\documentclass{amsart}
\usepackage{psfig}
\usepackage{epsf}
\usepackage{amssymb, amsmath,amscd}

\newcommand{\N}{\mathbb N}
\newcommand{\Q}{\mathbb Q}
\newcommand{\R}{\mathbb R}
\newcommand{\Z}{\mathbb Z}

\newcommand{\Ge}{\varepsilon}

\newcommand{\Gs}{\sigma}

\theoremstyle{plain}
\swapnumbers
\newtheorem{thm}{\itshape Theorem}[section]

\newtheorem{lem}[thm]{\itshape Lemma}
\newtheorem{olem}[thm]{\itshape Obvious Lemma}
\newtheorem{prop}[thm]{\itshape Proposition}

\newtheorem{cor}[thm]{\itshape Corollary}
\theoremstyle{definition}

\newtheorem{conj}[thm]{\itshape Conjecture}

\theoremstyle{remark}
\newtheorem{rem}[thm]{Remark}
\newcommand{\abz}{\newline\indent}

\newcommand{\cA}{\mathcal A}
\newcommand{\cD}{\mathcal D}
\newcommand{\cK}{\mathcal K}
\newcommand{\cP}{\mathcal P}
\newcommand{\lk}{\operatorname{lk}}
\newcommand{\rk}{\operatorname{rk}}
\newcommand{\sminus}{\smallsetminus}
\def\<{\langle}
\def\>{\rangle}

\begin{document}
\title{Finite Type Invariants of Classical and Virtual Knots}
\author[M. Goussarov]{Mikhail Goussarov}
\address{M. Goussarov, 
POMI, Fontanka 27,
St. Petersburg, 191011, Russia}
\email{goussar@pdmi.ras.ru}
\author[M. Polyak]{Michael Polyak}
\address{M.~Polyak, 
School of Mathematics, Tel-Aviv University,
69978 Tel-Aviv, Israel}
\email{polyak@math.tau.ac.il}
\author[O. Viro]{Oleg Viro}
\address{O.~Viro, Department of Mathematics, Uppsala University S-751 06
Uppsala, Sweden;\abz
POMI, Fontanka 27, St.Petersburg, 191011, Russia.}
\email{oleg@math.uu.se}
\begin{abstract}
We observe that any knot invariant extends to virtual knots.  The
isotopy classification problem for virtual knots is reduced to an
algebraic problem formulated in terms of an algebra of arrow diagrams.
We introduce a new notion of finite type invariant and show that the
restriction of any such invariant of degree $n$ to classical knots is
an invariant of degree $\le n$ in the classical sense. A universal
invariant of degree $\le n$ is defined via a Gauss diagram formula.
This machinery is used to obtain explicit formulas for invariants of
low degrees. The same technique is also used to prove that any finite
type invariant of classical knots is given by a Gauss diagram
formula. We introduce the notion of $n$-equivalence of Gauss
diagrams and announce virtual counter-parts of results concerning
classical $n$-equivalence.
\end{abstract}

\maketitle
\section{Virtualization}

Recently L. Kauffman introduced a notion of a virtual knot,
extending the knot theory in an unexpected direction. We show here
that this extension  motivates a new approach to finite type
invariants. This approach leads to new results both for virtual and classical
knots.

\subsection{Diagrams and Gauss Diagrams}
\begin{figure}[hbt]
\centerline{\epsffile{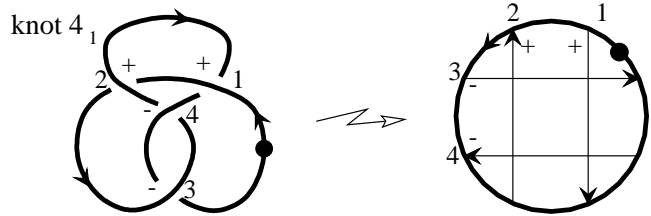}}
\caption{A diagram of the figure eight knot and its corresponding Gauss
diagram. }
\label{f2}
\end{figure}

Knots (smooth simple closed
curves in $\R^3$)  are usually presented by knot diagrams which are
generic immersions of the circle into the plane enhanced by information on
overpasses and underpasses at double points. A generic immersion of a
circle into the plane is characterized by its {\it Gauss diagram,\/}
which consists of
the circle together with the preimages of each double point of the
immersion  connected
by a chord. To incorporate the information on overpasses and
underpasses, the chords are oriented from the upper branch to the lower
one. Furthermore, each chord is equipped with the sign of the
corresponding double point (local writhe number). See Figure \ref{f2}.
The result is called a {\it Gauss diagram\/} of the knot.

A Gauss diagram is usually considered up to orientation preserving
homeomorphism of the underlying 
circle. The Gauss diagram defines (up to isotopy of $S^2$) a knot
diagram on the sphere, i.e., a knot diagram embedded into $S^2$ via the
embedding $\R^2\to S^2$. Given a knot diagram on the sphere, a knot
diagram on the plane can be recovered  modulo a finite ambiguity (this
involves the choice of which connected component of the diagram
contains the point at infinity), but the underlying  knot itself is
recovered uniquely up to isotopy.

Thus Gauss diagrams can be considered as an alternative way to
present knots. Of course, they cannot compete with knot diagrams in
creating a visual impression of knots, but Gauss diagrams are simpler
from the combinatorial point of view and provide numerous advantages
when we want to calculate knot invariants.

Unfortunately, not every picture which looks like a Gauss diagram is
indeed a Gauss diagram of some knot. Moreover, this is not easy to
recognize.
There is an obvious algorithm \cite{CE} for checking
this, which is just the result of attempting to draw the corresponding knot
diagram. However this requires a considerable amount of  effort.

\subsection{Virtual Knots}
The starting point for the present work is the idea that for some
purposes it is easier just to ignore the problem of whether a Gauss
diagram represents a knot, rather than trying to solve it.
This gives rise to a generalization of classical knot theory by
replacing true knots with objects which generalize Gauss diagrams of knots,
but which are not necessarily associated to a knot. Of course, these objects
are to be considered up to an appropriate equivalence, which imitates
knot isotopy.

Although we had been led to this generalization by the internal logic of our previous
research on combinatorial formulae for Vassiliev knot invariants, as
soon as we formulated it, we recognized that we had
rediscovered the theory of {\it
virtual knots,\/} which was announced last year by Louis Kauffman in
several talks \cite{Ka}. Our main contribution to this newborn theory is
to turn it into a useful tool for studying classical knots.

\begin{figure}[htb]
\centerline{\epsffile{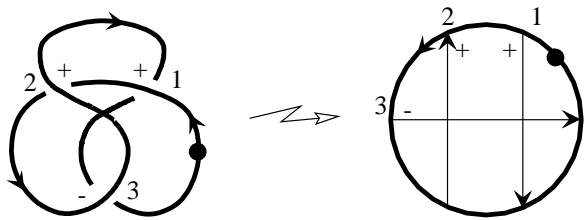}}
\caption{A diagram of a virtual knot with three real crossings and one virtual
crossing and its corresponding Gauss diagram.}
\label{virtf2}
\end{figure}

A {\it virtual knot diagram\/} is a generic
immersion of the circle into the plane, with double points divided into real
crossing points and virtual crossing points, with the real crossing points 
enhanced by information on overpasses and underpasses (as for
classical knot diagrams). At a virtual crossing the branches are not
divided into an overpass and an underpass. The {\it Gauss diagram\/} of a
virtual knot is constructed in the same way as for a classical
knot, but all virtual crossings are disregarded,  see Figure
\ref{virtf2}.

Virtuality of  virtual knots is manifest through the fact that while
the diagram
of a real knot is a picture describing a curve in $\R^3$, a virtual knot
diagram apparently does not pertain to any familiar 3-dimensional geometric object.
However we would like to keep speaking about virtual knots in the same
way we that speak about real knots: a virtual knot is that thing presented by a
virtual knot diagram. To resolve this ambiguity, we introduce moves
on virtual knot diagrams similar to the moves of real knot diagrams
which happen during an isotopy of a knot, and we will use the term virtual knot
to denote an equivalence class of virtual knot diagrams under these moves.
Two knot diagrams represent the same virtual knot, if one can be
obtained from the other by a sequence of these moves.
This agrees with the tradition of
classical knot theory, where the term knot is often taken to refer to
the isotopy class
of a knot.

\subsection{Reidemeister Moves and Virtual Moves} 
\begin{figure}[hbt]
\centerline{\epsffile{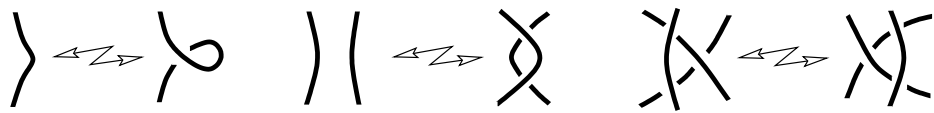}}
\caption{Reidemeister moves.}
\label{moves1}
\end{figure}

As is well-known, when a
knot changes by a generic isotopy, its diagram undergoes a sequence of {\it
Reidemeister moves\/} of one of the three types shown in Figure \ref{moves1}.

\begin{figure}[hbt]
\centerline{\epsffile{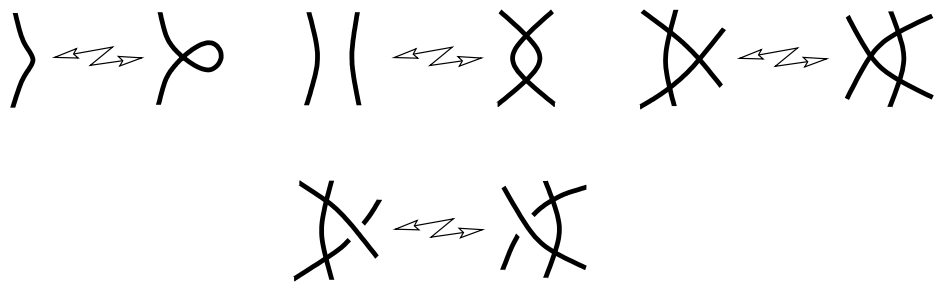}}
\caption{Virtual moves.}
\label{moves2}
\end{figure}

A diagram of a virtual knot can undergo the same Reidemeister moves, as
well as the moves shown in Figure \ref{moves2}. These additional moves are
called {\it virtual moves.\/}  The first three of them are 
versions of the Reidemeister moves, but with  virtual crossings in
place of crossings.
The last one looks like the third Reidemeister move, but involves two
virtual crossings and one usual crossing. 

\begin{figure}[hbt]
\centerline{\epsffile{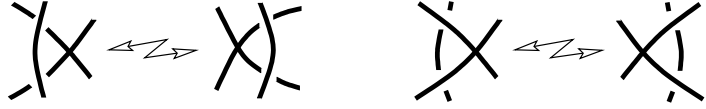}}
\caption{Forbidden moves.}
\label{forbmove}
\end{figure}

Similar moves, but with two real crossings and one virtual crossing
(shown in Figure \ref{forbmove}) are forbidden. If one allows these
moves, this makes the theory trivial: any virtual knot diagram can be
unknotted by a sequence of moves shown in Figures \ref{moves1},
\ref{moves2} and \ref{forbmove}, see Section \ref{sn-EVK.3}.

As mentioned above, a virtual knot is a class of virtual knot
diagrams consisting of diagrams which can be transformed into each other by sequences of
Reidemeister and virtual moves. A sequence of this kind is called a
{\it virtual isotopy.\/}

Virtual moves do not affect Gauss diagrams. On the other hand,
virtual moves allow one to move the interior of any arc which does not
pass through a real crossing quite arbitrarily. Therefore we obtain the

\begin{thm}
A Gauss diagram defines a virtual knot diagram  up to virtual moves.
\end{thm}

\begin{figure}[hbt]
\centerline{\psfig{file=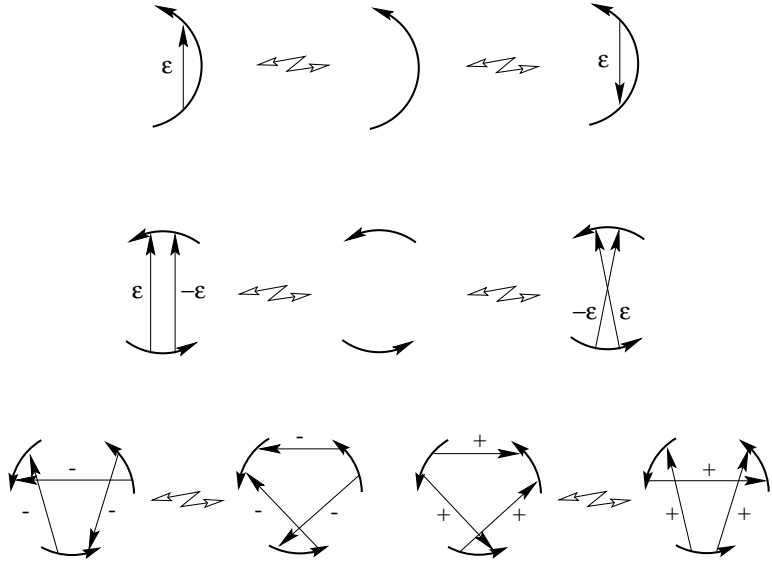,width=4in,silent=}}
\caption{Moves of Gauss diagrams corresponding to Reidemeister moves.}
\label{movesCD}
\end{figure}

This means that a virtual knot (modulo Reidemeister and virtual moves) is
equivalent to the corresponding Gauss diagram considered up to moves which are
the counter-parts of Reidemeister moves for Gauss diagrams, see Figure
\ref{movesCD}. Since in a Gauss diagram all the orientations and the
cyclic ordering of the endpoints of arrows are essential, each type of
Reidemeister moves splits. In Figure \ref{movesCD}, all moves
corresponding to the first and second Reidemeister moves are shown in
the top and middle rows, respectively. There are eight moves
corresponding to the third Reidemeister move, but we only show two of
them in the bottom row. As \"Ostlund \cite{Ost} showed, the remaining
six moves are unnecessary. That is, any sequence
of moves of a Gauss diagram can be replaced by a sequence of moves
appearing in Figure \ref{movesCD}. Although in \cite{Ost} this  is proved
for Gauss diagrams of knots, the same proof works for virtual knots. 

\subsection{Kauffman's Results on Extending Knot Invariants to Virtual
Knots}\label{sKR} Kauffman \cite{Ka}, \cite{Ka1} has proved that many knot
invariants extend to invariants of virtual knots. In particular, the
notions of knot
group, quandle and rack, and the bracket polynomial, all extend in a
straightforward manner. He also announced that with the use of a
"virtual framing" there are extensions of all quantum link invariants
and the large collection of their corresponding Vassiliev invariants.

The extensions are done in a formal way, disregarding the original
topological nature of these invariants. For example, the knot group, which is
defined for classical knots as the fundamental group of the knot
complement, is extended via a formal construction of a Wirtinger
presentation. This construction can be written down in terms of a
Gauss diagram as follows.

Let $G$ be a Gauss diagram. If we cut the circle at each arrowhead
(forgetting arrowtails),   the circle of $G$ is divided
into a set of arcs.  To each of these arcs there corresponds a generator
of the group. Each arrow gives rise to a relation. Suppose the sign of an
arrow is $\Ge$, its tail lies on an arc labelled $a$, its head is the
final point of an arc labelled $b$ and  the initial point of an arc
labelled $c$. Then we assign to this arrow the relation $c=a^{-\Ge}ba^{\Ge}$.
The resulting group is called the {\it  group of the Gauss diagram.\/}
One can easily check that it is invariant under the Reidemeister moves
shown in Figure \ref{movesCD}. Moreover, the group
system\footnote{Recall that the group system of a knot is the knot
group together with the class of subgroups which are conjugate to the
subgroup which is generated by a meridian and longitude.} also extends.
For the meridian,
take the generator corresponding to any of the arcs. To write down the
longitude, we go along the circle starting from this arc and write
$a^{\Ge}$, when passing the head of an arrow whose sign is $\Ge$ and
whose
tail lies on the arc labelled $a$.

The notion of quandle \cite{Jo} is extended in the same way
as the knot group: the generators remain the
same, but each group relation $c=a^{-\Ge}ba^{\Ge}$ is replaced with the
corresponding quandle relation 
$c=a\triangleright^{\Ge}b$.

\subsection{Knots Versus Virtual Knots}\label{sKVVK}
Any diagram of a classical knot can be considered to be a virtual knot
diagram. A virtual isotopy can turn it into a diagram with
virtual crossings, and then back again to a real knot diagram. Thus
virtual isotopy is
a new relation among classical knots, which apriori
could differ from classical isotopy. However this is not the case.

\begin{thm}[Virtual Isotopy Implies Isotopy]\label{VIII} {\mdseries
(See also Kauffman \cite{Ka}.)} Virtually
isotopic classical knots are isotopic.
\end{thm}

\begin{proof}
The group system extends to virtual knots. Hence it is preserved under
virtual isotopy, and virtually isotopic knots have isomorphic group
systems. Now recall that the group system is a complete knot invariant:
knots with isomorphic group systems are isotopic. \end{proof}

Any invariant of virtual knots is obviously an invariant of classical
knots. On the other hand, by Theorem \ref{VIII}, any invariant of
classical knots can be extended to an invariant of virtual knots.
Nevertheless, 
for some invariants it is not easy to choose a natural extension. Even
for the linking number, the extension to virtual 2-component links is not
unique (see Section \ref{sL} below). A similar situation occurs for the
degree 2 Vassiliev knot invariant considered in Section \ref{sGDF2}.

These examples are based on the same phenomenon. Unlike a classical
knot, a virtual knot cannot be turned upside down. 
A rotation of a classical knot by the angle $\pi$ around a horizontal
line reverses all arrows of its Gauss diagram, while  their signs do
not change. An application of this operation to a Gauss diagram of a
virtual knot gives rise to a virtual knot which may be non isotopic to
the original one. The composition of this operation with an invariant
of virtual knots may be another invariant. 

\begin{figure}
\centerline{\epsffile{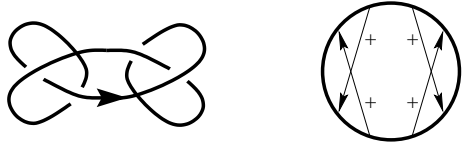}}
\caption{A virtual knot with different upper and lower groups.}
\label{uplowgr}
\end{figure}

A striking manifestation of this phenomenon comes from the knot group.
Instead of the upper Wirtinger presentation of a knot group, which was
generalized to virtual knots by Kauffman (see the preceding section),
let us use the lower Wirtinger presentation, i.e. compose Kauffman's
construction with arrows reversal. We will call these groups the
{\em upper\/} and  the {\em lower\/} virtual knot groups, respectively.
A virtual knot with different upper and lower groups is shown in Figure
\ref{uplowgr}. The upper group of this knot is isomorphic
to the group of the trefoil knot, while its lower group is $\Z$.

These examples may create an impression that virtual knot theory is
more cumbersome than the classical knot theory. However, this is not the
case. Due to its larger class of objects, the theory of virtual knots
provides more flexibility. This leads to significant simplification,
especially in the theory of finite type invariants.

\subsection{Long Knots}\label{sLK}
By a (classical) {\it long knot\/} we mean a smooth embedding $\R\to\R^3$
which coincides with the standard embedding outside a compact set.

An isotopy of long knots is a smooth isotopy in the class of embeddings
above.  In the classical knot theory, long knots are introduced for
purely technical reasons, since adding the point at infinity turns a
long knot into a knot in the sphere $S^3$ and this construction
establishes a one-to-one correspondence between the isotopy classes of
long knots and the isotopy classes of knots.

Given a diagram of a long knot, the corresponding Gauss diagram is the
line parameterizing the knot, together with arcs connecting the preimages of
each crossing.
As in the case of {\it closed\/} classical knots
considered above, the arcs are oriented from the upper branch to the
lower one and equipped with signs which are equal to the local writhe
numbers of the corresponding crossing points. Each oriented signed
arc is called an {\it arrow.\/} 

A {\it virtual long knot diagram\/} is a generic immersions $\R\to\R^2$
with double points divided into real and virtual
crossing points, where real crossing points are enhanced by information
on overpasses and underpasses, as in a classical knot diagram. At a
virtual crossing the branches are not divided into an overpass and an
underpass. The {\it Gauss diagram\/} of a virtual long knot is constructed
in the same way as for a classical long knot, but all the virtual
crossings are disregarded. A {\it virtual long knot\/} is a class of
 diagrams which can be transformed into each other  by sequences
of Reidemeister and virtual moves (shown in Figures \ref{moves1} and
\ref{moves2}).

\begin{figure} [htb]
\centerline{\epsffile{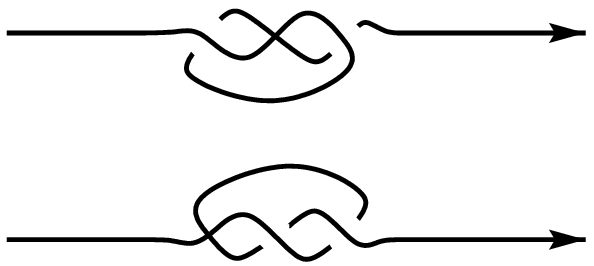}}
\caption{}
\label{fNILKs}
\end{figure}

Surprisingly, virtual long knots differ from virtual knots. That is,
there is no
one-to-one correspondence between the virtual isotopy classes of virtual long
knots and virtual knots. Addition of a point at infinity of
the plane of the diagram
turns a diagram of a virtual long knot into a virtual knot diagram
on $S^2$. Removal of a point from the complement of the diagram yields a
virtual knot diagram on $\R^2$. One can easily prove that 
its virtual isotopy class does not depend on the choice of this
point. Thus we have a
natural map from the set of virtual long knots to virtual knots. This map
is surjective but not injective.  The simplest pair of virtual long
knots which are not virtually isotopic, but give rise to isotopic
virtual closed knots is shown in Figure \ref{fNILKs}. These virtual
long knots are distinguished by the invariant $v_{2,2}$ defined in
Section \ref{sGDF2} below.
The upper diagram can be transformed into the lower diagram by moving
the underpassing arc of the leftmost crossing through the point at
infinity. If these were classical knots, these transformation could be
replaced by moving the same arc under the rest of the diagram by a
sequence of Reidemeister moves. In our case this is impossible, since
we cannot apply the move of Figure \ref{forbmove}.

As a result, the theories for long and closed virtual knots are quite
different, although their restrictions to usual knots coincide.

\subsection{Links}\label{sL} For links, the basic notions of the virtual
theory are introduced in  a straightforward way. The only change is
that the underlying circle of a Gauss diagram is replaced with several
circles. Simple examples show that in many respects  it is richer than
the classical one and sometimes looks surprising. For instance, for
2-component links there are two independent versions of the linking
number. The invariant $\lk_{1/2}$ may be computed as a sum of signs of
real crossings where the first component passes over the second one.
Similarly, one can define $\lk_{2/1}$ by exchanging the components in
the definition of $\lk_{1/2}$ above. 

String links are related to links like long knots are to knots.
A classical $n$-component {\em string link\/}
is a smooth embedding of a disjoint union of $n$ copies of $\R$ into
$\R^3$ which coincides with the standard embedding outside a compact
set.
Here, by the standard embedding, we mean the one given by the formula
$t\mapsto(t,k,0)$, with $t\in\R$  and  $k=1,\dots,n$.
All the basic notions of the virtual theory extend naturally to string
links. A Gauss diagram in this case consists of $n$ parallel lines
and signed arrows with end points on these lines.

\section{Finite Type Invariants}\label{sFOI}

\subsection{Crossing Virtualization Versus Crossing
Change}\label{sCVVCC} In the
realm of virtual knots there is an elementary operation which does not
exist for classical knots. A real crossing can be turned into a virtual
one. In terms of Gauss diagrams it looks even simpler: we erase an
arrow. This operation simplifies the knot in the sense that after applying 
it a sufficient number of times we eventually get to the unknot.

In the classical knot theory an operation with this property is widely
used. This is the crossing change. However, it is more complicated in several
ways. First, in order to turn a knot into the unknot, one must apply this
move according to a certain pattern, say making the diagram descending or
ascending, whereas virtualizing crossings leads to the unknot
automatically. Second, the unknotting by crossing changes involves a
choice: even if we have chosen to proceed towards a
descending diagram, we still have to choose a point at which the
descent
begins. The result considered as a diagram depends on this choice.
Unknotting by virtualization eliminates these technically unpleasant
problems. Third, crossing changes do not diminish the number of
crossings, while each virtualization diminishes the number of real
crossings. Finally,  virtualization is more elementary than
crossing changing, since a crossing change can be presented as the
composition of one crossing virtualization and the inverse of another.

A more general operation defined on Gauss diagrams of virtual
knots is  passage to subdiagrams. Here $D'$ is a subdiagram of $D$ if
all the arrows of $D'$ belong to $D$. In this case we write $D'\subset
D$.

\subsection{Classical Finite Type Invariants}\label{sCFTI}The standard
theory of finite type invariants is based on crossing change as the
basic modification.

Recall that a function $\nu$ defined on the set of knot isotopy types
and taking values in an abelian group $G$ is said to be a finite type
invariant of degree $\le n$, if for any  knot
diagram $D$ and $n+1$ crossing points $d_1$, $d_2$, \dots, $d_{n+1}$
of $D$
\begin{equation}\label{eqDEGn}\sum_{\Gs}(-1)^{|\Gs|}\nu(D_{\Gs})=0.
\end{equation}
Here $\Gs=\{\Gs_1,\dots,\Gs_{n+1}\}$ runs over $(n+1)$-tuples of zeros
and ones, $|\Gs|$ is the number of ones in $\Gs$, and $D_\Gs$ is the
diagram obtained from $D$ by switching all crossings $d_i$ with
$\Gs_i=1$. 

This description can be simplified by extending a knot invariant to knots
with double points (called also {\it singular knots\/}). A double point appears when at a crossing point one
moves the upper branch downwards through the lower
branch. The knot with a double point is identified with the formal
difference between the two knots obtained by resolving the double point
in two ways. This can be formulated as the following formal relation:
\begin{equation}\label{eqCROSS}
\vcenter{\epsffile{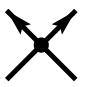}}=
\vcenter{\epsffile{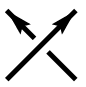}} -
\vcenter{\epsffile{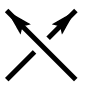}}
\end{equation}
Double points are depicted with thick points, so as to distinguish them from
virtual crossings.

Any knot invariant extends to formal linear combinations of knot
diagrams by
linearity. Under the identification in \eqref{eqCROSS}, the alternating sum in the left
hand side of equality \eqref{eqDEGn} becomes the value of
$\nu$ on a knot with $n+1$ double points. Thus a knot invariant has
degree at most $n$ if its extension vanishes on every singular knot
having at least
$n+1$ double points.

\subsection{A New Notion of Finite Type Invariant}\label{sNNFTI}
The counter-part in the virtual theory of the notion of finite type
invariant can be described as follows. We introduce a
new kind of crossing, which is called {\it semi-virtual.\/} At a
semi-virtual crossing there are still over- and under-passes. In a
diagram a semi-virtual crossing is shown as a real one, but surrounded
by a small circle. Semi-virtual crossings are related to the other types of
crossings  by the following formal relation:
\begin{equation}\label{eqSVCROSS}
\vcenter{\epsffile{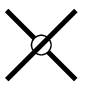}}=
\vcenter{\epsffile{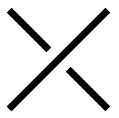}} -
\vcenter{\epsffile{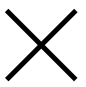}}
\end{equation}

In a Gauss diagram a semi-virtual crossing is presented by a dashed
arrow. The relation \eqref{eqSVCROSS} becomes
\begin{equation}\label{eq.semi-virt}
\vcenter{\epsffile{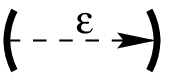}}\ =\
\vcenter{\epsffile{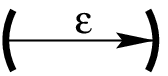}}\  -\ 
\vcenter{\epsffile{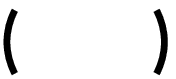}}
\end{equation}

Let $D$ be a virtual knot diagram and $\{d_1,\dots,d_n\}$ be an $n$-tuple
of its real crossings points.
For an $n$-tuple $\Gs=\{\Gs_1,\dots,\Gs_n\}$ of zeros and ones, define
$D_\Gs$ to be the diagram, obtained from $D$, by switching all the crossings
$d_i$, with $\Gs_i=1$, to virtual crossings.                           
Denote by $|\Gs|$ the number of ones in $\Gs$.
The formal alternating sum
$$\sum_\Gs (-1)^{|\Gs|} D_\Gs$$
is called a diagram with $n$ semi-virtual crossings.
We depict the corresponding alternating sum of Gauss diagrams by
the Gauss diagram of $D$ with 
all the arrows associated to $\{d_1,\dots,d_n\}$ being dashed.
This agrees with the convention \eqref{eq.semi-virt} on  semi-virtual
crossings. 

Denote by $\cK$ the set of virtual knots.
Let $\nu:\cK\to G$ be an invariant of virtual knots with values in an
abelian group $G$.
Extend it to $\Z[\cK]$ by linearity.
We say that $\nu$ is an {\em invariant of finite type}, if for some
$n\in\N$ it vanishes for any virtual knot $K$ with more than $n$
semi-virtual crossings.
The minimal such $n$ is called the {\em degree\/} of $\nu$.

Note that \eqref{eqSVCROSS} and \eqref{eqCROSS} imply
\begin{equation}\label{eqCR2SVCR}
\vcenter{\epsffile{ocross.eps}}=
\vcenter{\epsffile{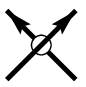}} -
\vcenter{\epsffile{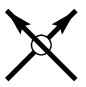}}
\end{equation}
It follows that for any finite type invariant of the
virtual theory, its restriction to classical knots
is a finite type invariant (of at most the same degree) in the classical
sense.

The definition of finite type invariants extends to virtual links in a
natural way.
A particularly simple example is given by the invariants $\lk_{1/2}$
and $\lk_{2/1}$ considered in Section \ref{sL}.
These invariants of 2-component virtual links have degree one.                                        

\subsection{The Algebra of Arrow Diagrams}\label{sFOI2}
An {\em arrow diagram\/} (on a circle) is an abstract diagram, which
consists of an oriented circle with pairs of distinct points connected
by dashed arrows.
Each arrow is equipped with a sign.
The {\em algebra of arrow diagrams\/} $\cA$ is the free abelian group
generated by all arrow diagrams. We call $\cA$ an algebra, because there is
indeed a natural multiplication in $\cA$ making it into an associative
algebra.\footnote{The product of arrow diagrams $A_1$, $A_2$ is the
sum (with appropriate multiplicities, cf. \cite{Ost}) of all
diagrams each of which is the union of subdiagrams isomorphic to $A_1$
and $A_2$.} The algebra $\cA\otimes\Q$ is isomorphic to the one
introduced in \cite{PV}. However in this paper we will not make use of
the multiplicative structure in $\cA$.

Denote the set of all Gauss diagrams by $\cD$.
Starting from any Gauss diagram we get an arrow diagram just by
making all its arrows dashed.
The extension of this map to $\Z[\cD]$ defines a natural isomorphism
$i:\Z[\cD]\to\cA$.

There is another important map $I:\cD\to\cA$, assigning to a Gauss
diagram $D$ the sum of all its subdiagrams and then making each of them
dashed:
$$I(D)=\sum_{D'\subset D}i(D').$$
\noindent 
Thus the map $I$ can described by the following symbolic formula:

\begin{equation}\label{eqI}
I: \vcenter{\epsffile{rarrow-ep.eps}}\ \mapsto
\vcenter{\epsffile{dasharrow-ep.eps}}\ +\
\vcenter{\epsffile{noarrow.eps}}.
\end{equation}

The reason for using the same dashed arrows both for semi-virtual
crossings in $\Z[\cD]$ and for arrows in $\cA$ becomes clear if one
compares formulas \eqref{eqI} and \eqref{eq.semi-virt}.

Extend  $I$ to $\Z[\cD]$ by linearity.

\begin{prop}\label{lem.iso} $I:\Z[\cD]\to\cA$ is an isomorphism.
The inverse map $I^{-1}:\cA\to\Z[\cD]$ is defined on
the generators of $\cA$ by  
$$I^{-1}(A)=\sum_{A'\subset A}(-1)^{|A-A'|}i^{-1}(A'),$$
\noindent where $|A-A'|$ is the number of arrows of $A$ which do not
belong to $A'$.  \qed
\end{prop}

A Gauss diagram is called
{\it semi-virtual\/} if each of its arrows is dashed.

\begin{cor}\label{basis}Semi-virtual diagrams form a basis of
$\Z[\cD]$.\qed
\end{cor}

\begin{rem}\label{rem}
 We can now explain an additional reason for the dual use of dashed arrows:
$I$ maps each semi-virtual Gauss diagram to the arrow diagram with the
same arrows.
This observation extends to diagrams containing both solid
and dashed arrows. Consider such a diagram $D$ as an
element of $\Z[\cD]$. Then each diagram appearing in $I(D)\in\cA$
contains all the dashed arrows of $D$.

Thus we see that $I$ can be interpreted  as  a presentation of a Gauss
diagram by a linear combination of semi-virtual diagrams.
\end{rem}

\subsection{The Polyak Algebra}\label{sFOI3}
The {\em Polyak algebra\/} is the quotient of $\cA$ 
by the following relations:
\begin{equation}\label{oneTprime}
\vcenter{\epsffile{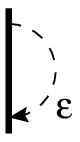}}\ =\ 0
\end{equation}
\begin{equation}\label{threeTprime}
\vcenter{\epsffile{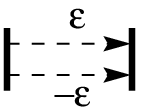}}\
+ \
\vcenter{\epsffile{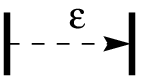}}\
+ \ \vcenter{\epsffile{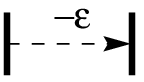}}\
= \ 0
\end{equation}
\begin{multline}\label{eightTprime}
\vcenter{\epsffile{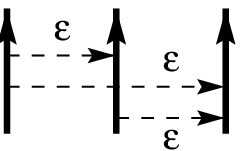}}\ + \ 
\vcenter{\epsffile{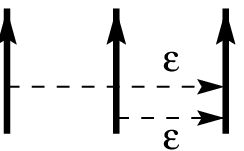}}\ +\ 
\vcenter{\epsffile{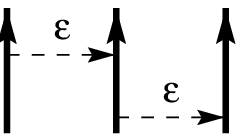}}\ +\ 
\vcenter{\epsffile{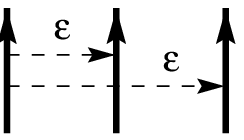}}\  = \\
\vcenter{\epsffile{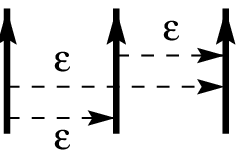}}\ + \ 
\vcenter{\epsffile{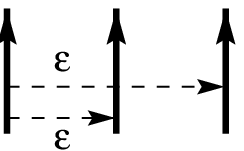}}\ +\ 
\vcenter{\epsffile{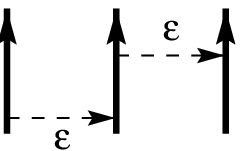}}\ +\ 
\vcenter{\epsffile{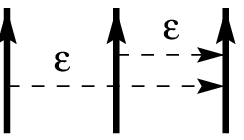}}\  
\end{multline}

Here we follow the common convention that the unshown parts of 
all diagrams involved in each of the relations coincide.
The embeddings of the shown parts into the whole diagrams should
preserve the orientations  in \eqref{eightTprime}.

The same relations define analogous algebras for long knots, links and
string links.   

The quotient of $\cA$ by the relations \eqref{oneTprime} -
\eqref{eightTprime} is an algebra, since the relations generate an
ideal of $\cA$, but we shall not go into further detail on this point.
Denote this algebra by
$\cP$. This algebra is closely related to the algebra $\cA$ introduced
in \cite{P}.

The isomorphism $I$ induces an isomorphism $I:\Z[\cK]\to\cP$ of
quotient algebras.
Indeed, the equivalence relation induced in $\cD$
by the Reidemeister moves shown in Figure \ref{moves1},
can be rewritten in terms of diagrams with semi-virtual crossings as
follows
\begin{equation}\label{oneT}
\vcenter{\epsffile{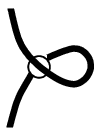}}\ =\ 0
\end{equation}
\begin{equation}\label{threeT}
\vcenter{\epsffile{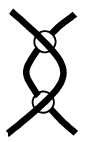}}\ +\ \vcenter{\epsffile{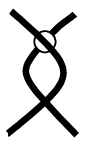}}
\ + \ 
\vcenter{\epsffile{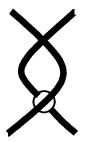}}\ =\ 0
\end{equation}
\begin{multline}\label{eightT}
\vcenter{\epsffile{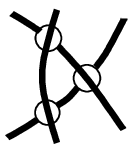}}\ +
\vcenter{\epsffile{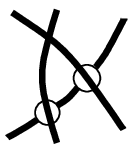}}\ +
\vcenter{\epsffile{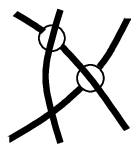}}\ +
\vcenter{\epsffile{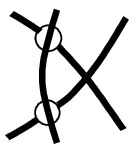}}\ = \\
\vcenter{\epsffile{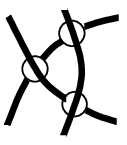}}\ +
\vcenter{\epsffile{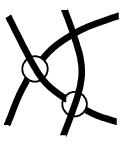}}\ +
\vcenter{\epsffile{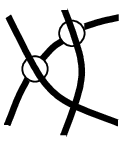}}\ +
\vcenter{\epsffile{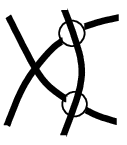}}
\end{multline}
\noindent
Note that the map $I$ turns \eqref{oneT} -- \eqref{eightT} into
\eqref{oneTprime} -- \eqref{eightTprime}.
Thus for any Gauss diagram $D$ of a virtual knot $K$, $I(D)$ defines
a $\cP$-valued invariant of $K$. 
Moreover, since the Gauss diagram determines $K$, the invariant $I(D)$
distinguishes virtual knots. Thus we obtain the following theorem.

\begin{thm}\label{thm.inv}                               
Let $D$ be any diagram of a virtual knot $K$.
The formula $K\mapsto I(D)\in\cP$ defines a complete invariant of virtual
knots.
\end{thm}                  

\subsection{The Truncated Algebras $\cP_n$ and the Universal Finite Type
Invariant}
\label{sFOI4}
Define the {\em truncated algebra\/} $\cP_n$ by putting $A=0$ for
any diagram $A\in\cP$ with more than $n$ arrows.
Denote by $I_n:\cK\to\cP_n$ the composition of $I:\cK\to\cP$ with the
projection $\cP\to\cP_n$.

Let $P$ be an abelian group and $p:\cK\to P$ be
a $P$-valued invariant of virtual knots.
We call $p$ a {\em universal invariant of degree $n$}, if for every
abelian group $G$ and every invariant $\nu:\cK\to G$ of degree at most
$n$ factors through $p$, i.e.\ there exists a map $\pi:P\to G$,
such that $\nu=\pi\circ p$.

\begin{thm}\label{thm.invn}
The map $I_n:\cK\to\cP_n$ defines a universal invariant of degree $n$.
\end{thm}
\begin{proof} Since, by Remark \ref{rem}, $I$ preserves all dashed
arrows, $I_n$ maps a knot with more than $n$ semi-virtual crossings
to zero. Therefore $I_n$ is an invariant of degree at most $n$.

Let $\nu:\cK\to G$ be an invariant of degree at most $n$.
We have to prove that $\nu\circ I^{-1}:\cP\to G$ factors through
$\cP_n$. Observe that, by Remark \ref{rem}, for any arrow diagram $A$
its image $I^{-1}(A)\in\Z[\cD]$ can be identified with the same diagram
$A$, but considered as a semi-virtual Gauss diagram. Since the degree
of $\nu$ is at most $n$, $\nu$ vanishes on each diagram with more than $n$
semi-virtual crossings. Therefore $\nu\circ I^{-1}$ vanishes on each
arrow diagram with more than $n$ arrows and $\nu\circ I^{-1}$ factors
through $\cP_n$. \end{proof}                                      

\begin{cor}
The space of $\Q$-valued invariants of degree at most $n$ is
finite-dimen\-s\-i\-o\-n\-al, of dimension equal to $\rk(\cP_n)$. 
It can be identified with the dual space $\cP_n^*$ of $\Q$-valued
linear functions on $\cP_n$.
\end{cor}

\section{Gauss Diagram Formulas for Finite Type Invariants}\label{sGDF}

\subsection{Gauss Diagram Formulas}\label{sGDF1}
Since the algebra $\cA$ has a distinguished basis, consisting of arrow
diagrams,  there is a natural orthonormal scalar product
$(\cdot,\cdot)$ on $\cA$.
Namely, on the generators of $\cA$ we put $(D_1,D_2)$ to be $1$, if
$D_1=D_2$, and $0$ otherwise and then extend $(\cdot,\cdot)$ bilinearly.
This allows us to define the pairing
$\<\cdot,\cdot\>:\cA\times\cD\to\Z$ in
the following way.
For any $D\in\cD$ and $A\in\cA$ put                                  
\begin{equation}\label{eq.bracket}
\<A,D\>=(A,I(D))=(A,\sum_{D'\subset D}i(D')).
\end{equation}
Informally speaking, we count subdiagrams of $D$ with weights, where
the weight of a diagram $D'$ is the coefficient of $i(D')$ in $A$.
 
In the case of Gauss diagrams corresponding to usual knots, this
pairing (in a slightly different form) was introduced in \cite{PV}
as a tool for writing down Gauss diagram formulas 
for knot invariants.
We will use it below for the same purpose in the framework of virtual
knots.

Using equation \eqref{eq.bracket} and Theorem \ref{thm.inv} it is easy
to see
that any $\Z$-valued invariant of finite type of virtual
knots can be obtained by a Gauss diagram formula
$\<A,\cdot\>:\cK\to\Z$ for some $A\in\cA$.
The maximal number of arrows of the diagrams in the linear combination
giving  $A$ is an upper bound for the
degree.
However, in general the expression $\<A,D\>$ depends on the choice
of the Gauss diagram $D$ of a virtual knot.
In our earlier work \cite{PV} we did not present any systematic
method for producing arrow polynomials $A$ which give invariants, and
we posed the following question: ``Which arrow polynomials define knot
invariants...?''
 We can now answer this question in the framework of the virtual theory: 
$A$ defines an invariant of degree at most $n$ if and only if
all diagrams in $A$ have at most $n$ arrows and $A$
satisfies the equations $(A,R)=0$, where $R$ runs over the left hand sides
of the relations $R=0$ defining $\cP_n$ in $\cA$, see Section
\ref{sFOI3}.

The general method for producing all invariants of degree $n$
requires a computation of the algebra $\cP_n$.
Some simple observations allow one to reduce this computation.
First, by a repeated use of \eqref{threeT}, one can eliminate arrows with
the negative sign. 
Second, since all the diagrams with more than $n$ arrows
vanish, for diagrams with $n$ arrows equations \eqref{threeT}
and \eqref{eightT} become simpler and contain only 2 and 6 terms
respectively, all of them with exactly $n$ arrows.
The simplified version of \eqref{threeT} implies
the following  rule for elimination of negative arrows in
diagrams with exactly $n$ arrows: 
if two such diagrams
differ only by the signs of $k$ arrows,
they differ in $\cP_n$ by multiplication by $(-1)^k$.
This allows one to drop the signs of arrows in diagrams with $n$
arrows, using the convention that a diagram with $k$ negative signs
of arrows is counted with coefficient $(-1)^k$.
The simplified version of \eqref{eightT} involves only diagrams with
exactly $n$ arrows and looks as follows:
\begin{multline}\label{sixT}%
\vcenter{\epsffile{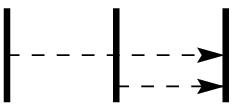}}\ +\ 
\vcenter{\epsffile{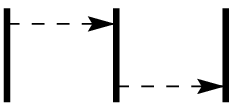}}\ +\ 
\vcenter{\epsffile{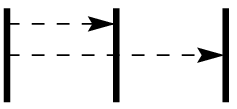}}\  = \\
\vcenter{\epsffile{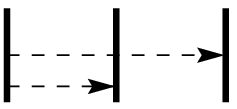}}\ +\ 
\vcenter{\epsffile{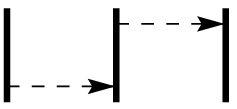}}\ +\ 
\vcenter{\epsffile{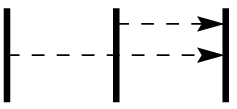}}\
\end{multline}
This $6T$-equation (introduced earlier in \cite{P}) is an oriented
version of the well-known $4T$-relation for chord diagrams.
The $4T$-relation can be recovered from \eqref{sixT} by repeated
use of the following formula which follows from \eqref{eqCR2SVCR}:
\begin{equation}\label{eqSVDASH}
\vcenter{\epsffile{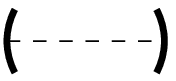}}\  =  \
\vcenter{\epsffile{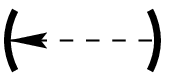}}\ + \
\vcenter{\epsffile{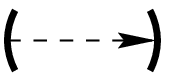}}. 
\end{equation}

\subsection{Invariants of Small Degree}\label{sGDF2}
Computation of $\cP_n$ for small $n$ can be done by hand and leads
to some interesting results.                                   
Similarly to the case of classical knots, there are no invariants
of degree one.
More surprisingly, the algebra $\cP_2$ is also trivial, so
there are no invariants of degree two! 
However, for long knots the corresponding algebra is 2-dimensional,
so there are two independent invariants $v_{2,1}$ and $v_{2,2}$ of degree 2.
These invariants are given by
\begin{equation}\label{eqv2}
v_{2,1}(\cdot)=\left\<\vcenter{\epsffile{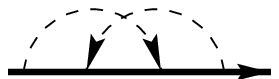}}\ ,\ \cdot\ \right\>;
\qquad v_{2,2}(\cdot)
=\left\<\vcenter{\epsffile{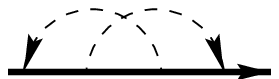}}\ ,\ \cdot\ \right\>
\end{equation}  
This illustrates a curious feature of the theory of virtual knots which
was discussed in Section \ref{sLK}.
For classical knots, there is one-to-one correspondence between the
isotopy classes of knots and long knots, hence any invariant of long
knots is an invariant of closed knots. We now see that for virtual
knots this is no longer true.

Another interesting feature of this theory is that many invariants,
which coincide for usual knots (due to the existence of certain symmetries),
are different on the larger class of virtual knots.
The invariants $v_{2,1}$ and $v_{2,2}$ provide a good illustration.

In degree three there is only one invariant, given by 
$$\left\<\vcenter{\psfig{figure=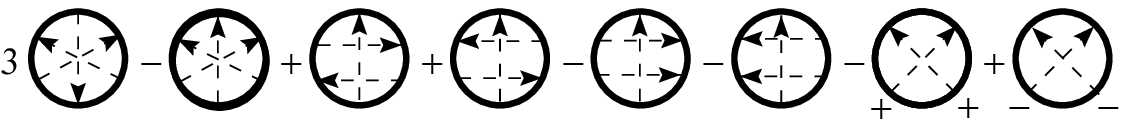,height=1cm,silent=}},\cdot \right\>.$$
It vanishes on real knots. 
Similarly to degree two, for long virtual knots there are several invariants
of degree three, which give the same degree three invariant of real knots. 
Here is an example of such an invariant:
$$\vcenter{\psfig{figure=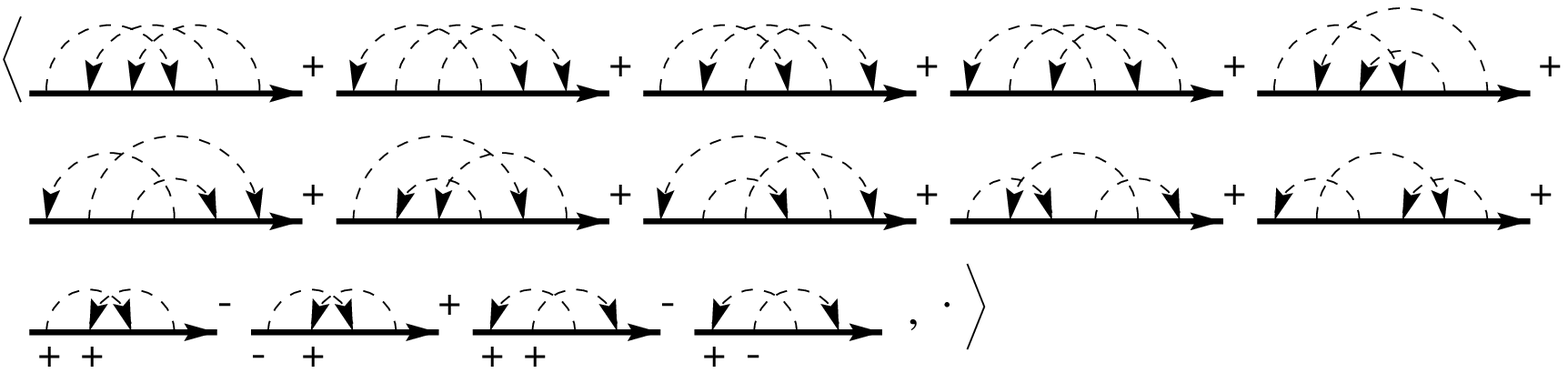,width=12.6cm,silent=}}$$ 

\subsection{The Case of Classical Knots}\label{sCCK}
Our work in this direction started with a search for combinatorial
formulas for finite type invariants of classical knots. The first
results were summarized in the paper \cite{PV} of the second and third
authors. There a class of combinatorial formulas similar to 
\eqref{eq.bracket} was introduced, and numerous special formulas of
this sort were found. In \cite{PV} we posed the following question:
``Can any Vassiliev invariant be calculated as a function of arrow
polynomials evaluated on the knot diagram?'' In the terminology used
above, an
arrow polynomial evaluated on the knot diagram is an expression of the
type given in \eqref{eq.bracket}.

This question has been answered in the affirmative by the first author. The
result is formulated as follows.

\begin{thm}[Goussarov]\label{GT}
Let $G$ be an abelian group and let $\nu$ be a $G$-valued 
invariant of degree $n$ of long (real) knots. Then there exists a
function $\pi:\cA\to G$  such that $\nu=\pi\circ I$
and such that $\pi$ vanishes on any arrow diagram with more than $n$ arrows.
\end{thm}

\begin{cor}\label{cGT}
Any integer-valued finite type invariant
of degree $n$ of long knots can be presented as $\<A,\cdot\>$, where
$A$ is a linear combination of arrow diagrams on a line with at most
$n$ arrows.
\end{cor}

The next section is devoted to the proof of this theorem.

To a large extent, the present paper was motivated by an analysis of
the proof of \ref{GT}, which originally was rather cumbersome.
The main difficulties in this proof were caused by the necessity of
requiring
all the numerous Gauss diagrams involved to be realizable.
The desire to get rid of this restriction motivated our interest in
virtual knots.
Indeed, for virtual knots the problem stated in \cite{PV} is solved
by Theorem \ref{thm.invn} above. The universal invariant of Theorem
\ref{thm.invn} is essentially $\sum_{A\in\cP_n}\<A,D\>A$, so any
$G$-valued invariant can be presented by a Gauss diagram formula.

Unfortunately, for classical knots the new technique
does not give a universal invariant. However it gives powerful and
simple machinery to generate Gauss diagram formulas for any invariant
which can be extended to a finite-type invariant of virtual knots. Hoping
for the best, we conjecture

\begin{conj}\label{conj}
Every finite-type invariant of classical knots can be extended to a
finite-type invariant of long virtual knots.
\end{conj}

This may require the  consideration of virtual framing.
The extension given by Kauffman \cite{Ka} \cite{Ka1} of numerous invariants
to virtual knots strongly supports this conjecture.

The main open problem concerning finite type invariants of classical
knots is whether such invariants distinguish  non-isotopic
knots. The positive solution of this problem would follow from the
positive solution of the corresponding problem for virtual knots.
By Theorems
\ref{thm.inv} and \ref{thm.invn}, the
latter can be reformulated in purely algebraic terms as the question
whether the natural map
$$\cP\to\lim_{\longleftarrow}\cP_n$$
is injective.

\section{Proof of Goussarov's Theorem}\label{sPGT}

\subsection{Scheme of the Proof}\label{sSchPrf}
The standard method  for calculation of a Vassiliev invariant
$\nu$ goes as follows (see \cite{V}). One picks a set of singular knots
which span (using relation \eqref{eqCROSS}) the free abelian group
generated by all nonsingular  knots. The invariant $\nu$ is determined
by its actuality table, i.e.\ its values  on this set. Given a knot
diagram, one unknots it, making it descending by a sequence of crossing
changes.   Under each crossing change the invariant jumps. The jump is
equal to the value of $\nu$ on the knot with a double point by
\eqref{eqCROSS}. Then each of these singular knots is deformed to a
knot with a single double point from the actuality table by an isotopy
and a sequence of crossing changes. The jumps of $\nu$ correspond to
knots with two double points. They are again deformed to the knots from
the actuality table. The process eventually stops when the number of
double points exceeds the degree of $\nu$ (by definition of the
degree).

In the proof of \ref{GT}, both the actuality table and the procedure of expansion
described above are made canonical. This is done by generalizing the
notion of a descending diagram to singular knots.\footnote{A similar
notion called almost monotone diagram was considered by
Bar-Natan \cite{B-N}, but the procedure of expansion in \cite{B-N}
involves some choices.}

More importantly, this is done in terms of Gauss
diagrams, so that the notion of descending diagram and the procedure of
expansion extend to virtual knots.   For a real
descending knot the isotopy class and hence the value of $\nu$
is determined by the part of the Gauss diagram encoding the double
points. For a virtual descending diagram the isotopy type is not
determined by this part of the Gauss diagram. Nevertheless we extend
$\nu$ to virtual descending diagrams literally in the same way. We
do not know whether the result is an invariant of virtual
knots.\footnote{A
positive answer to this question would imply Conjecture \ref{conj}.}
However, for our purposes, this formal extension turns out to be
sufficient.

Next we use the isomorphism $I^{-1}:\cA\to\Z[\cD]$ (see Proposition \ref{lem.iso})
to define $\pi:\cA \to G$ as $\nu\circ R$.
Some special properties of the extended map $\nu:\Z[D]\to G$ are then
used to prove that $\pi$ vanishes on diagrams with more than $n$
chords.

\subsection{Descending Singular Diagrams}\label{sDSD}
On a diagram of a long virtual singular knot each double point  is naturally
equipped with a sign. Indeed, the branches at a double point are
ordered and the sign is the intersection number of the branches
(taken in this order). On a Gauss diagram of a long singular knot, 
each double point is shown by a dashed chord equipped with the above sign.
A diagram $D'$ is called a {\em subdiagram\/} of a diagram $D$ if $D'$
consists of all the chords and some arrows of $D$.

Recall that a diagram of a real long knot is descending if going along
the knot in the positive direction we pass each crossing first going over and
then under. In terms of Gauss diagrams it means that all the arrows
are directed to the right.

\begin{figure}[htb]
\centerline{\psfig{figure=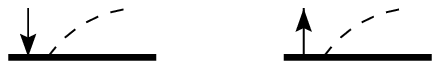,height=.8cm,silent=}}
\caption{}
\label{forbsit}
\end{figure}

We now extend this notion to virtual long knots with double points. We
still require that all the arrows are directed to the right. There is
also an additional condition:  there is no chord whose left endpoint
has  an endpoint of an arrow as immediate left neighbor. In other
words, the situations shown in Figure \ref{forbsit} are forbidden.

\begin{figure}[htb]
\centerline{\epsffile{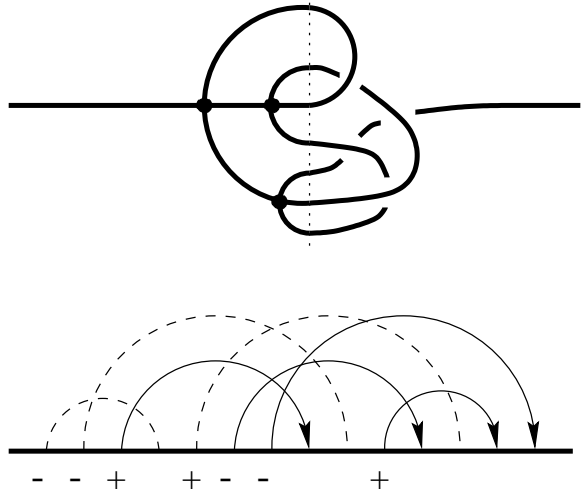}}
\caption{A descending long real knot and its Gauss diagram.}
\label{fDESC}
\end{figure}

A real long knot with a Gauss
diagram of this type can be presented by a diagram such that
\begin{itemize}
\item all the double points are in the left half-plane,
\item all the crossings are in the right half-plane,
\item the intersection of the diagram with the left
half-plane is an embedded tree,
\item the intersection with the right half-plane is an ordered collection of
arcs; each of them is descending and lies below all the
previous ones.
\end{itemize} See Figure \ref{fDESC}.

\begin{lem}\label{lem.C}
Let $D_1$ and $D_2$ be Gauss diagrams of real descending long knots and
let $\nu$ be an invariant of long knots. If the chord parts of $D_1$ and
$D_2$ coincide then $\nu(D_1)=\nu(D_2)$.
\end{lem} 
\begin{proof}
One can see that the isotopy class of a real descending long knot is
determined by the chord part of its Gauss diagram. Indeed, the chord
part  determines the tree in the left half-plane
(recall that the chords have signs, which define the embedding
locally). The rule for connecting the endpoints of the tree by arcs in
the right half-plane is determined by the mutual position of the chords
in the Gauss diagram. Since the diagram is descending, the connection
by the arcs is unique up to isotopy. 
\end{proof} 
  
\subsection{Reduction to Descending Diagrams}\label{sRDD}
There is an algorithm for expressing the  Gauss diagram of a long knot with
double points as a linear combination of descending diagrams. This
algorithm
consists of steps of two types. At each step, one inspects the Gauss
diagram from the left to the right looking for the first fragment where the
diagram fails to be descending. Such a fragment may either be a bad
arrow or a bad chord. An arrow is bad if it is directed to the left.
A chord is bad if an immediate left neighbor of its left
endpoint is an endpoint of an arrow, as in Figure \ref{forbsit}.

In the case of a bad arrow the step of the algorithm is the replacement of
the diagram with the sum of two diagrams according to the  formula
$$\vcenter{\epsffile{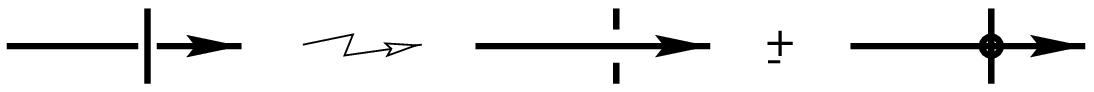}}$$
In terms of Gauss diagrams this replacement is as follows:
$$\vcenter{\epsffile{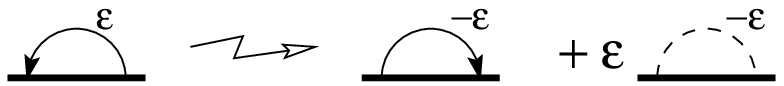}}$$

\begin{figure}[hptb]
\centerline{\psfig{figure=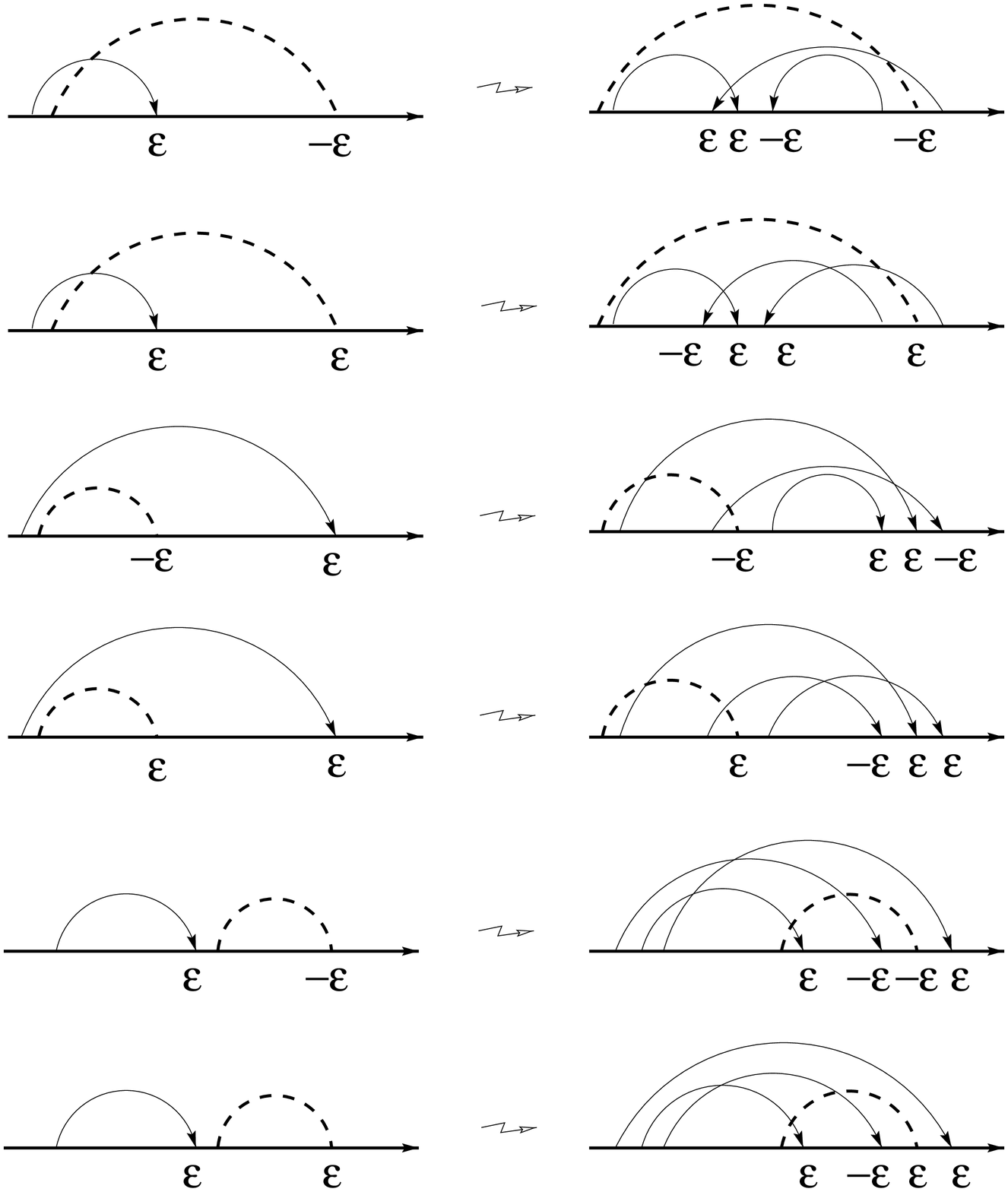,width=9cm,silent=}}
\caption{}
\label{pullback}
\end{figure}

In the case of a bad chord the step of the algorithm is the pulling of the
crossing over or under the appropriate branch by isotopy: 
$$\vcenter{\epsffile{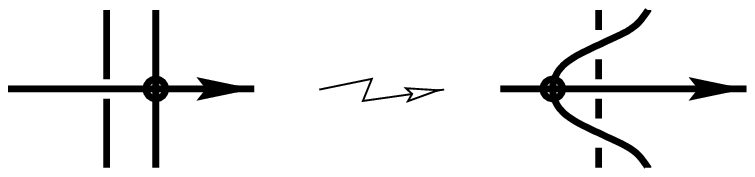}}$$
$$\vcenter{\epsffile{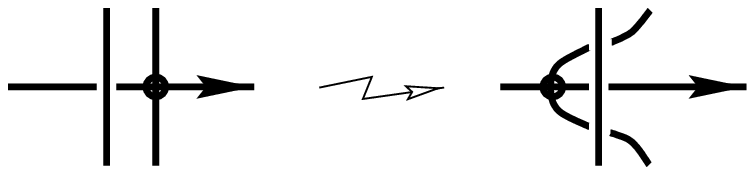}}$$
In terms of Gauss diagrams, this corresponds to one of the
transformations shown in Figure \ref{pullback}. The different cases in
Figure \ref{pullback} correspond to different  orientations and
possible orderings  of the three arcs.

Since we deal with an invariant of degree $n$, the diagrams with more
than $n$ chords are disregarded. Thus, when one applies a
step of the algorithm to a bad
arrow in a diagram with $n$ chords, the summand with $n+1$ chords
disappears. Denote by $\cD_n$ the free abelian group
generated by Gauss diagrams of virtual long singular knots with at most
$n$ chords (note that $\Z[\cD]=\cD_0\subset\cD_n$).
We will think of a step of our algorithm as of an operator acting on
$\cD_n$. Denote this operator by $P$. By the definition of $P$, for any
descending Gauss diagram $D$ we have $P(D)=D$.

\begin{lem}\label{lem.B1} For any diagram $D\in\cD_n$ there exists $m$ such
that $P^m(D)$ is a sum of descending diagrams.
\end{lem}

\begin{proof}
Let $l(D)$ be the number  of chords of $D$ which have one of the
endpoints
to the left of the first bad fragment. As is easy to see,
$l(D')\ge l(D)$ for each diagram in the expansion of $P(D)$.
However the number of such chords 
in a non-descending diagram is at most $n$. Therefore it
suffices to prove that the diagram cannot change infinitely many times
in subsequent iterations of $P$ without changing $l$.

Consider the number of arrowheads on the ray to the left of the left endpoint
of the $(l(D)+1)$-th  chord. For any diagram involved in
the expansion of $P(D)$  this number is not greater than that for $D$. If it
is the same for one of these diagrams, then it has less arrowtails on
the same ray. This can happen only finitely many times.
\end{proof}

\subsection{The Extension of $\nu$ and the Construction of $\pi$}\label{sEtV}Denote by 
$\cD_n^{re}$  the subgroup  of  $\cD_n$ generated by Gauss diagrams of
{\em real\/} long singular knots. Any finite type invariant of
classical knots of degree at most $n$ extends to $\cD_n^{re}$ by
linearity.

\begin{olem}\label{lem.A} The operator $P:\cD_n\to\cD_n$ preserves
$\cD_n^{re}$.  The restriction of $P$ to $\cD_n^{re}$ preserves any
invariant of degree at most $n$.\qed
\end{olem}

We now extend an invariant $\nu$ of degree at most $n$ to virtual
descending diagrams. For any such diagram $D$  there exists a
descending diagram $D^{re}$ of a real knot with the same double points, i.e.,
the same chord part of the Gauss diagram. $D^{re}$ can be obtained by turning
all the virtual crossings of $D$ into appropriate real ones. Put
$\nu(D)=\nu(D^{re})$. By Lemma \ref{lem.C}, $\nu(D^{re})$ does not
depend on the choice of $D^{re}$.

Next we extend $\nu$ to all virtual diagrams. By Lemma \ref{lem.B1}, for any
diagram $D\in\cD_n$ there exists $m$ such that $P^m(D)$ is a sum of
descending diagrams and hence $\nu(P^m(D))$ is already defined. Put
$\nu(D)=\nu(P^m(D))$. Lemma \ref{lem.A} implies that on $D^{re}$ this
agrees with the initial definition of $\nu$. Since $P^{m+1}(D)=P^m(D)$
we get:

\begin{olem}\label{P preserves nu} The operator $P:\cD_n\to\cD_n$ preserves
$\nu$, i.e.\ $\nu\circ P=\nu$.  \qed
\end{olem}

 We are now in a position to construct the map $\pi:\cA\to G$
of Theorem \ref{GT}. Define $\pi:\cA\to G$ as the
composition
$$\begin{CD}\cA@>I^{-1}>>\Z[\cD]\subset \cD_n@>\nu>>G,
\end{CD}$$
where $I^{-1}$ is the isomorphism of Proposition \ref{lem.iso}
and $\nu$ is the extension of the
original finite type invariant to $\cD_n$. Then for any diagram $D$ of
a long knot
$$\nu(D)=\pi(I(D))=\sum_{D'\subset D}\pi(i(D')).$$

In order to prove Theorem \ref{GT}, we must show that
$\pi(A)=0$ for any arrow diagram $A$ with more than $n$ arrows.
The rest of this section is devoted to the proof of this fact .

\subsection{The Analogues of $\cD_n$ and $P$ for Arrow Diagrams}\label{sTAD}
The algebra $\cA$ of arrow diagrams on the line
is generated by diagrams consisting of the line and dashed arrows
(oriented signed arcs). Consider now diagrams which in 
addition to arrows also contain dashed chords, i.e.\ unoriented signed
arcs.
Denote by $\cA_n$ the free abelian group generated by such diagrams
with at most $n$ chords. 

The maps $i,I:\Z[\cD]\to\cA$ defined in Section \ref{sFOI2} on Gauss
diagrams without chords extend to isomorphisms $i,I:\cD_n\to\cA_n$.
The chord parts of the diagrams remain intact under both $i$ and $I$,
while the arrows are dealt with as in Section \ref{sFOI2}.

We now define an operator $Q:\cA_n\to\cA_n$, which is an analogue
of $P$.  A diagram $A\in\cA_n$ is called
descending, if $i^{-1}(A)$ is descending. Put $Q(A)=A$ if $A$ is
descending. Otherwise, find the leftmost bad fragment of $A$ (the
notion of a bad fragment is borrowed from $\cD_n$ via $i$). If it is a
bad arrow, we define $Q(A)=iPi^{-1}(A)$. If it is a bad chord, put
$Q(A)=\sum A'$ where the sum runs over all the subdiagrams of
$iPi^{-1}(A)$, each of which contains all the arrows not shown in Figure
\ref{pullback}, all the chords and at least one more arrow. In other
words, we sum up all seven subdiagrams of $iPi^{-1}(A)$ which
contain all the arrows and chords also belonging to $A$ plus at least
one more arrow.

\begin{rem}\label{remQ}Observe that in both cases, we sum up all  subdiagrams of diagrams
in $iPi^{-1}(A)$ which are not subdiagrams of $A$, but contain all 
arrows of $A$ except for the arrow involved into the bad fragment.
The arrows of $A$ which are not in the leftmost bad fragment play a
passive role in the construction of $Q$: if $A'$ is a
subdiagram of $A$ obtined by removing arrows which are not in the
leftmost bad fragment, then $Q(A')$ is obtained from $Q(A)$ by removing
the same arrows from each of the summands.
\end{rem}

\begin{olem}\label{Q increases arr+ch} For any diagram $A\in\cA_n$, the
total number of arrows and chords in each diagram appearing in $Q(A)$
is at least the total number of arrows and chords in $A$.\qed
\end{olem}

\begin{lem}\label{lem.B2} For any diagram $A\in\cA_n$, there exists $m$ such
that $Q^m(A)$ is a sum of descending diagrams.
\end{lem}

The proof of this Lemma is completely analogous to the proof of Lemma
\ref{lem.B1}.\qed

\begin{lem}\label{lem.E}
For any non-descending diagram $D\in\cD_n$, there is a splitting $I(D)=U+V$ with
$U,V\in\cA_n$ such that 
\begin{equation}\label{I(P(D))}
I(P(D)) = Q(U)+V
\end{equation}
and such that $U=i(D)+U'$, where $U'$ is a sum of diagrams each of which has
fewer arrows than $D$.
\end{lem}

\begin{proof} Let $U$ be the sum of all the subdiagrams of $i(D)$
which include the first bad fragment of $i(D)$. These
subdiagrams contain the same bad fragment as the whole diagram $i(D)$.
As follows from Remark \ref{remQ}, $Q(U)$ is the sum of all subdiagrams
of diagrams in $iP(D)$ which are not subdiagrams of $i(D)$.
Then $V$ is the sum of
the subdiagrams of $i(D)$ which do not contain the arrow from the bad
fragment (in the case of a bad chord, this is the arrow shown on the left
hand side of Figure \ref{pullback}) and these subdiagrams of $i(D)$
remain  unchanged, when one applies $P$ to $D$. Thus $I(P(D)) = Q(U)+V$.
\end{proof}

\begin{lem}\label{Q preserves pi} The operator $Q:\cA_n\to\cA_n$ preserves
$\pi$, i.e.\ $\pi\circ Q=\pi$.
\end{lem}

\begin{proof} Let $A\in\cA_n$ be a diagram and $D=i^{-1}(A)$.
Let us prove that $\pi(Q(A))=\pi(A)$  by induction on the number of
arrows in $A$. If this number equals $0$, then $A$ is descending and
$Q(A)=A$ by definition of $Q$.
Suppose inductively that the statement is correct for any diagram whose number of
arrows is less then the number of arrows in $A$ and let us prove the statement
for $A$. Apply $\pi$ to \eqref{I(P(D))}:
$$\pi\circ Q(U)+\pi(V)= \pi\circ I\circ P(D)=\nu\circ P(D)$$
By Lemma \ref{P preserves nu} and the definition of $\pi$
$$\nu\circ P(D)=\nu(D)=\pi\circ I(D)=\pi(U)+\pi(V).$$
Thus $\pi\circ Q(U)=\pi(U)$.
By the induction assumption, $\pi\circ Q(U')=\pi(U')$, where $U'=U-A$
(as in Lemma \ref{lem.E}), and  we obtain
the desired equality $\pi(Q(A))=\pi(A)$. This completes the
induction step.\end{proof}

\begin{lem}\label{lem.-1} Let $A\in\cA_n$ be a descending diagram such
that the total number of arrows and chords in $A$ is greater than $n$.
Then $\pi(A)=0$.
\end{lem}
\begin{proof} Let $D=i^{-1}(A)$.
By the definition of $\pi$ and Proposition \ref{lem.iso},
$$\pi(A)=\nu\circ R(A)=\sum_{D'\subset D}(-1)^{|D-D'|}\nu(D').$$
Since any subdiagram $D'$ of $D$ is descending and has
the same chord part, $\nu(D')=\nu(D)$ by the construction of $\nu$.
Therefore
$$\pi(A)=\left(\sum_{D'\subset D}(-1)^{|D-D'|}\right)\nu(D).$$
As one can easily check by induction on the number of arrows in $A$,
the sum in parentheses is equal to 1 if $A$ has no arrows and is 0
otherwise. Since all the diagrams in $\cA_n$ have at most $n$
chords and the total number of arrows and chords in $A$ is greater
than $n$,  it has at least one arrow. Hence $\pi(A)=0$.
\end{proof}

\begin{lem}\label{lem.-0} Let $A\in\cA_n$ be a  diagram such
that the total number of arrows and chords in $A$ is greater than $n$.
Then $\pi(A)=0$.
\end{lem}
\begin{proof}
Let $m$ be the number which exists for $A$ by Lemma
\ref{lem.B2}. By Lemma \ref{Q preserves pi}, $\pi(A)=\pi(Q^m(A))$.  By 
Lemma \ref{Q increases arr+ch}, the expansion of $Q^m(A)$ contains only
descending diagrams with the total number of chords and arrows greater
than $n$. Then by Lemma \ref{lem.-1}, $\pi(A)=0$. \end{proof}

This concludes the proof of Theorem \ref{GT}.\qed

\section{$n$-Equivalence}\label{sn-EVK}                                                                 

\subsection{$n$-Trivial Gauss Diagrams}\label{sn-EVK.1} Let $D$ be a
Gauss diagram and let its arrows be colored with $n$ colors. Consider
all subdiagrams which can be obtained from $D$ be removing all
arrows colored with one or several colors. If all arrows of each of
these diagrams can be removed by the second Reidemeister moves, then
the coloring is said to be destroying.

A Gauss diagram, based on a union of several disjoint segments, is
called {\it $n$-trivial\/} if it admits a destroying coloring with $n+1$
colors.

The property of $n$-triviality does not change if one reverses 
the orientation of some of the segments. It also does not change if one simultaneously
reverses the orientations or the signs of all  arrows
connecting two segments.

\subsection{$n$-Variations}\label{sn-EVK.2}  On a Gauss diagram $D$,
choose several segments which do not contain an endpoint of any arrow.
Adjoin to $D$ arrows of an $(n-1)$-trivial  Gauss diagram based on the
chosen segments. This transformation of $D$ is called an {\it
$n$-variation.\/}

It is easy to see that an addition of any number of arrows is a $1$-variation.
On a virtual knot diagram an addition of an arrow can be realized as
follows: $$\vcenter{\psfig{figure=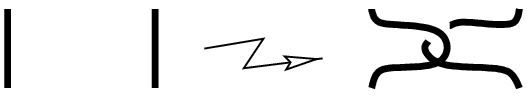,height=.7cm,silent=}}.$$

\begin{figure}[htb]
\centerline{\epsffile{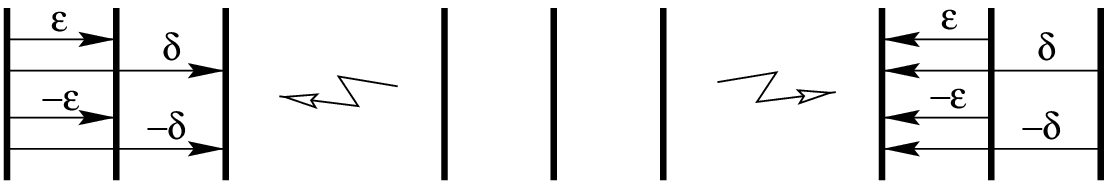}}
\caption{}
\label{2-var}
\end{figure}

In Figure \ref{2-var} we show the simplest 2-variations. To get the
corresponding destroying coloring, one colors all arrows connecting the
first two strings with one color and the other arrows with the other
color. It is easy to see that these 2-variations do not change
$\lk_{i/j}$.  

On a virtual knot diagram these 2-variations can be realized as 
$$\vcenter{\psfig{figure=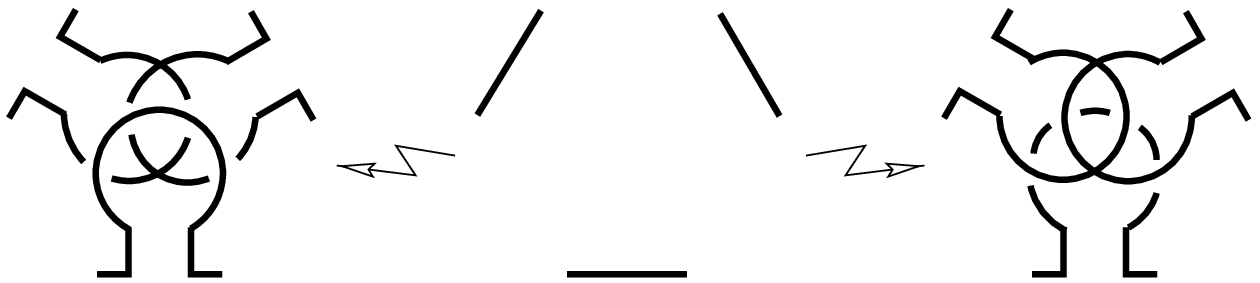,height=1in,silent=}}$$
Observe that these modifications coincide, up
to isotopy, with the forbidden moves of Figure \ref{forbmove}. 

Some obvious properties of $n$-variations are: \begin{enumerate}
\item An $n$-variation is a $k$-variation for any $k<n$.
\item Composition of several $n$-variations is an $n$-variation. 
\end{enumerate}

A less obvious property is: the result  of  an isotopy followed by an
$n$-variation can be presented as the result of other $n$-variations 
followed by an isotopy.

The following proposition is a key property of $n$-variations. 
\begin{prop}\label{}After  any $n$-variation, one can apply another
$n$-variation such that the final result is the initial diagram, up to
a sequence of  second Reidemeister moves.
\end{prop}

\subsection{$n$-Equivalence}\label{sn-EVK.3}  Two Gauss diagrams are said
to be {\it $n$-equivalent\/} if they can be transformed to each other
by a sequence of isotopies and $(n+1)$-variations. For example, the
Gauss diagrams shown in Figure \ref{2-var} are 1-equivalent. Moreover,
one can prove that any two 1-equivalent Gauss diagrams can be
transformed to each other by a sequence of isotopies and the
2-variations of Figure \ref{2-var}. A 1-equivalence class of string
links is completely determined by the invariants $\lk_{i/j}$. Therefore
any two closed virtual knots are 1-equivalent and can be transformed to
each other by a sequence of the Reidemeister moves and the forbidden moves of
Figure \ref{forbmove}.

The transition from the use of Gauss diagrams to that of 
$n$-equivalence classes  yields better results when the set of Gauss
diagrams is equipped with a natural multiplication. For example, in the
case of virtual string links (and, in particular, long knots)
$n$-equivalence classes form a group.

In the cases of virtual knots and (closed) links the set of
$n$-equivalence classes has a more complicated algebraic structure.
As in the case of classical links, the following trick works. Consider
a virtual string link with $2n$ strings. It can be turned into a closed
one by adding $n$ arcs from above and below (this generalizes the plat
presentation of a link with braids replaced by  string links).
This gives rise to a map from the group of $n$-equivalence classes of
virtual string links to the set of $n$-equivalence classes of closed
virtual links. This map is a double coset factorization. From the left
we quotient out by string links which become $n$-equivalent to the trivial
one by adding only arcs from below, and  from the right, similarly with 
arcs from above. Both sets of string links give rise to subgroups which
are not normal in general.

The value of a finite type invariant of degree $\le n$ depends only on
the $n$-equivalence class. Usually, in a non-group situation, invariants of
degree $\le n$ do not separate all the $n$-equivalence classes. For
instance, virtual closed knots do not admit an invariant of degree 2.

\end{document}